\documentclass[11pt,leqno]{amsart}

\usepackage{amssymb,amsmath,amsthm,amsfonts}
\usepackage{latexsym}

\newcommand{\Hmm}[1]{\leavevmode{\marginpar{\tiny%
$\hbox to 0mm{\hspace*{-0.5mm}$\leftarrow$\hss}%
\vcenter{\vrule depth 0.1mm height 0.1mm width \the\marginparwidth}%
\hbox to 0mm{\hss$\rightarrow$\hspace*{-0.5mm}}$\\\relax\raggedright #1}}}
\newcommand{\nc}{\newcommand}
\nc{\les}{\lesssim}
\nc{\nit}{\noindent}
\nc{\nn}{\nonumber}
\nc{\D}{\partial}
\nc{\diff}[2]{\frac{d #1}{d #2}}
\nc{\diffn}[3]{\frac{d^{#3} #1}{d {#2}^{#3}}}
\nc{\pdiff}[2]{\frac{\partial #1}{\partial #2}}
\nc{\pdiffn}[3]{\frac{\partial^{#3} #1}{\partial{#2}^{#3}}}
\nc{\abs}[1] {\lvert #1 \rvert}
\nc{\cAc}{{\cal A}_c}
\nc{\cE}{{\cal E}}
\nc{\cF}{{\mathcal F}}
\nc{\cP}{{\cal P}}
\nc{\cV}{{\cal V}}
\nc{\cQ}{{\cal Q}}
\nc{\cGin}{{\cal G}_{\rm in}}
\nc{\cGout}{{\cal G}_{\rm out}}
\nc{\cO}{{\cal O}}
\nc{\Lav}{{\cal L}_{\rm av}}
\nc{\cL}{{\cal L}}
\nc{\cB}{{\cal B}}
\nc{\cZ}{{\cal Z}}
\nc{\cR}{{\cal R}}
\nc{\cT}{{\cal T}}
\nc{\cY}{{\cal Y}}
\nc{\cX}{{\cal X}}
\nc{\cXT}{{{\cal X}(T)}}
\nc{\cBT}{{{\cal B}(T)}}
\nc{\vD}{{\vec \mathcal{D}}}
\nc{\efield}{\mathcal{E}}
\nc{\vE}{{\vec \efield}}
\nc{\vB}{{\vec \mathcal{B}}}
\nc{\vH}{{\vec \mathcal{H}}}
\nc{\ty}{{\tilde y}}
\nc{\tu}{{\tilde u}}
\nc{\tV}{{\tilde V}}
\nc{\Pc}{{\bf P_c}}
\nc{\bx}{{\bf x}}
\nc{\bX}{{\bf X}}
\nc{\bXYZ}{{\bf XYZ}}
\nc{\bY}{{\bf Y}}
\nc{\bF}{{\bf F}}
\nc{\bS}{{\bf S}}
\nc{\dV}{{\delta V}}
\nc{\dE}{{\delta E}}
\nc{\TT}{{\Theta}}
\nc{\dPsi}{{\delta\Psi}}
\nc{\order}{{\cal O}}
\nc{\Rout}{R_{\rm out}}
\nc{\eplus}{e_+}
\nc{\eminus}{e_-}
\nc{\epm}{e_\pm}
\nc{\eps}{\varepsilon}
\nc{\vnabla}{{\vec\nabla}}
\nc{\G}{\Gamma}
\nc{\w}{\omega}
\nc{\mh}{h}
\nc{\mg}{g}
\nc{\vphi}{\varphi}
\nc{\tlambda}{\tilde\lambda}
\nc{\be}{\begin{equation}}
\nc{\ee}{\end{equation}}
\nc{\ba}{\begin{eqnarray}}
\nc{\ea}{\end{eqnarray}}

\nc{\g}{\gamma}
\nc{\ol}{\overline}

\newtheorem{theo}{Theorem}[section]
\newtheorem{defin}{Definition}[section]
\newtheorem{prop}{Proposition}[section]
\newtheorem{lem}{Lemma}[section]

\newtheorem{rmk}{Remark}[section]

\nc{\pT}{\partial_T}
\nc{\pz}{\partial_z}
\nc{\pt}{\partial_t}
\nc{\la}{\langle}
\nc{\ra}{\rangle}
\nc{\infint}{\int_{-\infty}^{\infty}}
\nc{\halfwidth}{6.5cm}
\nc{\figwidth}{10cm}

\nc{\nlayers}{L} \nc{\nsectors}{M}
\nc{\indicator}{\mathbf{1}}
\nc{\Rhole}{R_{\rm hole}}
\nc{\Rring}{R_{\rm ring}}
\nc{\neff}{n_{\rm eff}}
\nc{\Frem}{F_{\rm rem}}
\nc{\Real}{\mathbb R}
\nc{\Z}{\mathbb Z}
\nc{\DD}{\Delta}
\nc{\cD}{\mathcal D}
\nc{\lnorm}{\left\|}
\nc{\rnorm}{\right\|}
\nc{\rnormp}{\right\|_{\ell^{p,\eps}}}
\nc{\rar}{\rightarrow}
\date{\today}

\begin{document}

\title[Near-linear dynamics in KdV]{Near-linear dynamics in KdV with periodic
boundary conditions}

\author{M.~B.~Erdo\smash{\u{g}}an, N.~Tzirakis, and V.~Zharnitsky}
\thanks{The authors were partially supported by NSF grants DMS-0600101 (B.~E.), DMS-0901222 (N.~T.), and DMS-0807897  (V.~Z.)}

\address{Department of Mathematics \\
University of Illinois \\
Urbana, IL 61801, U.S.A.}

\email{berdogan@uiuc.edu \\ tzirakis@math.uiuc.edu\\ vzh@uiuc.edu}

\maketitle

\begin{abstract}
Near linear evolution in  Korteweg de Vries (KdV) equation with periodic boundary conditions is
established under the assumption of high frequency initial data. This result is obtained by the method of normal form reduction.
\end{abstract}

\section{Introduction}
This articles investigates the behavior of a class of solutions with high frequency initial data of Korteweg de Vries (KdV) equation,
\[
v_t =  6 v v_x -  v_{xxx},
\]
with periodic boundary conditions $v(x+2\pi)=v(x)$. We show that, see Theorem~\ref{theo:main} below, these solutions evolve
near linearly, ({\em i.e.} like solutions of $v_t =  -  v_{xxx}$) for large times.

On the real line, near linear behavior  for dispersive PDEs, such as nonlinear Schr\" odinger equation, nonlinear Klein-Gordon equations, KdV, {\em etc.}, could be expected. Indeed, in that case high frequency solutions will disperse over a large subset of the real line weakening the nonlinearity.
For example, under some conditions, one can extend the $L^1-L^{\infty}$ dispersive estimates for the linear Schr\"odinger equation to NLS, see {\em e.g.} \cite{Naum}.
In the focusing case, linear evolution could be destroyed by the mass concentration phenomenon as it leads to larger nonlinear effects. However, such concentration cannot occur in the case of mass subcritical nonlinearity. In short, there are two major reasons why on the real line, the evolution of high frequency solutions in the mass subcritical NLS case should be near linear: {\em dispersive decay} and {\em absence of collapse}.

For the KdV on the torus or a circle (periodic boundary conditions), the linear solution is periodic in space and time and, thus, 
one does not have dispersive decay. It is also generally believed that the solutions of KdV on the torus will not be approximated by the linear solutions as time goes to infinity.
 Therefore, it is somewhat surprising that, as we show in this paper, the evolution is still near linear on a finite but large time scale.
One can argue that a hint towards this behavior comes from Bourgain's discovery of Strichartz estimates for periodic case \cite{Bou1, Bou2}.
In some way, the effect of dispersion for the periodic problem can be interpreted as averaging of the nonlinearity over high frequencies.

On the torus there are other reasons, such as resonances, which could prevent linear behavior. For NLS, see \cite{erdogan_vz}, such  resonances cause faster phase rotation while the behavior is still linear.

Our results are also motivated by the scattering problems for dispersive PDEs. On the real line there are many results on scattering, which show that  nonlinear solutions tend to the linear ones as time goes to infinity. On the torus, however, one does not expect scattering. For example, the absence of scattering  was proved rigorously for the cubic NLS on the two dimensional torus in \cite{ckstt}. Our statement is different since we only claim linear behavior for large but finite time scale for a special class of high frequency solutions. On the other hand, our near linear solutions provide some scattering like behavior.

Although KdV with periodic boundary conditions is completely integrable, our methods do not rely on integrability. We only use the conservation of momentum, energy, and Hamiltonian. An interesting question is whether integrability 
structure can be used to obtain more precise results on near linear evolution and on a larger time scale.

Our work also suggests a new mechanism of formation of  the so-called rogue waves.
Rogue waves (also called freak and giant waves) correspond to large-amplitude waves appearing
on the sea surface ``from nowhere''. In the scientific literature, the following amplitude criterion for the rogue wave
is usually used: its height should exceed the significant wave height (on the sea surface) by about a factor of two
\cite{KharifPelin}.

There is  a vast literature on rogue waves, see {\em e.g.} the survey paper \cite{KharifPelin} and
references therein, and many explanations have been proposed. Some scenarios  involve
\begin{itemize}
\item probabilistic approach -- rogue waves are considered as rare events in the framework of Rayleigh statistics
\item linear mechanism -- dispersion enhancement (spatio-temporal focusing)
\item nonlinear mechanisms -- in approximating models ({\rm e.g.} NLS or KdV), for some special
initial data large amplitude waves can be created.
\end{itemize}

Linear mechanism of rogue wave formation is simpler since there are various solutions leading to large amplitudes, 
while nonlinear mechanism requires very special initial data. On the other hand, linear approximations are valid  in the small
amplitude limit which is restrictive.
This article shows that for KdV the linear and nonlinear mechanisms can be combined into one since
we describe a large subset of initial data for which the solutions of KdV equation behave near linearly.

Regarding the boundary conditions, our choice of periodic boundary conditions is not the most
realistic one but appropriate for a model problem. Indeed, while the sea surface is not periodic,
one observes more or less similar pattern over large areas.

We finally mention that it would be best to observe near-linear dynamics for
the full water wave problem, however, it is a considerably harder problem which
will be addressed in future work. We also limit our study to the one dimensional problem.

\section{Main Results}
We consider KdV equation
\be
v_t =  6 v v_x -  v_{xxx},
\ee
with periodic boundary conditions $v(x+2\pi)=v(x)$ and we assume $v \in H^1(S^1)$.
In this case, KdV is well-posed \cite{Temam} and can be written in Hamiltonian form
\[
v_t = \frac{d}{dx} \frac{\D H}{\D v},
\]
where the Hamiltonian is given by
\be
H(v) = \int_{-\pi}^{\pi} \left (\frac{1}{2} v_x^2 + v^3 \right ) dx
\label{eq:hamph}
\ee
and $\frac{\D H}{\D v}$ denotes $L^2-$gradient of $H$, representing the Fr\' echet derivative of $H$
with respect to the standard scalar product on $L^2$. 
We  also need to consider linear part of KdV
\[
v_t + v_{xxx}=0,
\]
with the solution given by
\[
v(x,t) = e^{Lt}v(x,0),
\]
where $L = -\D_{xxx}$.

While KdV possesses infinitely many conserved quantities, we  use
the first  three: the above Hamiltonian, linear momentum
\be
P = \int_{-\pi}^{\pi} v(x) dx
\ee
and kinetic energy
\be
K = \int_{-\pi}^{\pi} v^2(x) dx.
\label{eq:kinetic}
\ee

\begin{theo}
Assume without loss of generality\footnote{One can reduce the case $P\neq 0$ to the zero momentum case $P=0$ by a simple transformation. } that $P=0$ and
\be
\|v(\cdot,0)\|_{H^1} \leq C \eps ^{-1},\,\,\,  \|v(\cdot,0)\|_{H^{-1/2}} \leq C \eps ^{1/2}.
\label{eq:thinspec}
\ee
for some $C>0$ and for sufficiently small $\eps>0$. Then for any $t \les \eps^{-\frac12 +}$
\[
\|v(\cdot,t) - e^{Lt} v(\cdot,0)\|_{L^2} \les \langle t \rangle \,\, \eps^{\frac12 -},
\]
where the implicit constant depends only on $C$ but not on $\eps$.
\label{theo:main}
\end{theo}

This Theorem follows from Theorem~\ref{theo:main_u} below, which is proved by applying near-identical
canonical transformations, so that the new Hamiltonian
flow is close to the linear one. This implies that the original Hamiltonian flow is also close
to the linear one.

\begin{rmk}\label{rmk:uni-bound}
Note that since the Hamiltonian \eqref{eq:hamph} and the kinetic energy \eqref{eq:kinetic} are conserved quantities, the bounds \eqref{eq:thinspec}
imply $|H(v(t))| \les \eps^{-2}$. This immediately implies a uniform bound in time $\|v\|_{H^1} \les \epsilon^{-1}$.
\end{rmk}

To prove our theorem we first apply the following transformation \cite{KappelerPoschel},
which is a weighted modification of Fourier transform
\be
v(x)=\sum_{n\neq 0}\sqrt{|n|}e^{inx}u(n),
\label{eq:fourier}
\ee
where $n\in \Z\backslash \{0\}$ and $u(n)$ is a bi-infinite sequence of complex numbers.
Since, $v(x)$ is real,
\[
u(-n)= \overline{u(n)}.
\]

In these new variables the Hamiltonian takes the form\footnote{Below we will omit
 the absolute value sign $|*|$ under the square root. It will be implicitly implied for
the rest of the paper. }
\begin{align}
H & = i \sum_{n>0} n^3 u(n) u(-n) + i \sum_{n_1+n_2+n_3=0} \sqrt{n_1 n_2 n_3} \,\,
u(n_1)u(n_2)u(n_3) \nn\\
&=: \Lambda_2 + H_3,
\label{eq:ham2}
\end{align}
where $\Lambda_2$ and $H_3$ are the quadratic and cubic parts of the Hamiltonian.
Equivalently, in order to deal with the summation over all $n\neq 0$, we can write
\[
\Lambda_2 = \frac{i}{2} \sum_{n\neq 0} n^3\sigma(n)u(n)u(-n),
\]
where $\sigma(n):=\text{sgn}(n)$.
In this formulation $u(m)$ and $u(-m)$, with $m=1, 2, ...$ are conjugated canonical
variables with the standard symplectic structure, so that
the Hamiltonian equations take the usual form
\[
\frac{du(m)}{dt} = \frac{\D H}{\D u(-m)}
\]
\[
\frac{du(-m)}{dt} = -\frac{\D H}{\D u(m)},
\]
where $m>0$. We also write these equations in a more compact form
\[
\frac{du(m)}{dt} = \sigma(m) \, \frac{\D H}{\D u(-m)},\,\, {\rm where} \,\, m\neq 0.
\]

It is straightforward to verify that these are the correct equations, by applying the change of variable (\ref{eq:fourier})
directly to KdV.

Now, we introduce a subset of $l^2$
\[
X_{\eps}^{\rho}=\left \{
u \in l^2: u(0) = 0, u(-n)= \bar u(n), \|u\|_{l^2} \leq \rho \sqrt{\eps}, \|u\|_{l^2_{3/2}} \leq \frac{\rho}{\eps}
\right \},
\]
where
$$
\|u\|_{\ell^2_s}^2=\sum_k  |k|^{2s} \, |u(k)|^2.
$$
We will also need the norm
$$
\|u\|_{\ell^p_s}^p=\sum_k  |k|^{ps} \, |u(k)|^p.
$$

Note that the hypothesis of Theorem~\ref{theo:main} is equivalent to $u\in X_{\eps}^{\rho}$ initially in time for some $\rho>0$.
 By Remark \ref{rmk:uni-bound}, for any $t>0$, $\|u(\cdot,t)\|_{l^2_{3/2}}\les \eps^{-1}$.
For the initial data in this subset we prove that the evolution is near linear.
\begin{theo}\label{theo:main_u}
Let $\rho > 0$ be fixed. Assume $u(\cdot,0) \in X_{\eps}^{\rho}$ for sufficiently small $\eps >0$.
Then for any $t \les \eps^{-\frac12 +}$, $u(\cdot,t)\in X_{\eps}^{2\rho}$ and
\be
\|u(n,t) - e^{in^3 t}u(n,0)\|_{l^2_{1/2}(n)} \les \langle t \rangle \,\, \eps^{\frac12-}. \label{eq:t_est}
\ee
\end{theo}
Theorem~\ref{theo:main} immediately follows from this one by applying the relation $u(n)=\hat v(n)/\sqrt{|n|}$.
To prove Theorem~\ref{theo:main_u}, we apply two canonical transformations $\Phi_{F_1}^1$, $\Phi_{F_2}^1 $,
see the next section, so that $u=u(q)=\Phi_{F_1}^1 \circ \Phi_{F_2}^1 (q)$. The new Hamiltonian is given by
\ba \nn
H(q) = H(u(q)) = \Lambda_2(q) + R(q),
\ea
where $R$ stands for the reminder terms, and the equations take the form
\be
\dot q(n) = i n^3 q(n) + E(q)(n),
\label{eq:qande}
\ee
where
\be \label{eq:Eqn}
E(q)(n)= \frac{\D }{\D q(-n)} R, \,\,\,\,\,n>0.
\ee
The transformation is near-identical in the following sense:
\begin{prop}
If $u \in X_{\eps}^{\rho}$ or $q \in X_{\eps}^{\rho}$, then
\be
\|u(q)-q\|_{l^2_s} \les \eps^{1-s},
\ee
where $s\in [0,3/2]$ and the implicit constant depends on $\rho, s$ but not on $\eps$.
In particular, for sufficiently small $\eps$, if $q\in X_{\eps}^{\rho}$, then
$u(q) \in X_{\eps}^{2\rho}$ and vice versa.
\label{prop:transf}
\end{prop}
The estimate for the error term is given by
\begin{prop}
If $q\in X_{\eps}^{\rho}$ then the error term satisfies
\be
\|E(q)\|_{l^2_s} \les \eps^{1-s-},
\ee
where $s\in [0,1/2]$ and the implicit constant depends on $\rho, s$ but not on $\eps$.

\label{prop:error}
\end{prop}
\begin{proof}[ Proof  of Theorem \ref{theo:main_u}]
The proof follows easily from Propositions \ref{prop:transf} and  \ref{prop:error}. Indeed,
multiplying (\ref{eq:qande}) with the integrating factor $e^{-in^3 t}$ and integrating from 0 to $t$, we obtain
\begin{align}
q(n,t)e^{-i n^3 t} - q(n,0) = \int_0^t e^{-i n^3 \tau} E(q)(n) d\tau.
\label{eq:qest}
\end{align}
Next, by taking the $\ell^2_s$ norm after multiplying both sides with $e^{in^3 t}$, we obtain
\[
\|q(n,t) - e^{in^3 t}q(n,0)\|_{l^2_s} =
\left\| \int_0^t e^{i n^3 (t-\tau)} E(q)(n) d\tau
\right\|_{l^2_s} \leq |t| \,\, \|E\|_{l^2_s} \lesssim |t| \,\, \eps^{1-s -}
\]
for $s \in [0,1/2]$ and $t \les \eps^{-\frac12 +}$.

Then, using the triangle inequality, we estimate, for $s\in [0,1/2],$
\begin{align*}
 &\|u(n,t) - e^{in^3 t}u(n,0)\|_{l^2_s}  \leq  \nn \\
&\leq \|u(n,t)-q(n,t) \| _{l^2_s} + \|q(n,t) - e^{in^3 t}q(n,0)\|_{l^2_s} + \nn \\ &
\|e^{in^3 t}q(n,0) -e^{in^3 t}u(n,0)\|_{l^2_s}  \,\, \lesssim \,\, \langle t \rangle \eps^{1-s -}.
\end{align*}
The first and the third estimates follow from Proposition \ref{prop:transf} while the second  follows from the estimate on the equation 
\eqref{eq:qest}.

This inequality for $s=0$ implies that $\|u(\cdot,t)\|_{\ell^2}\leq 2\rho \sqrt{\eps}$ for $t \les \eps^{-1/2+}$, while the conservation of Hamiltonian implies that $\|u(\cdot,t)\|_{\ell^2_{3/2}}\leq 2\rho  \eps^{-1}$. Therefore,   $u$ stays in $X^{2\rho}_{\eps}$ up to the time $t \les \eps^{-1/2+}$. This is important since our estimates for the canonical transformations are only valid in the ball $X^{C\rho}_\eps$.

Moreover, for $s=1/2$, the last inequality   gives \eqref{eq:t_est}.
This ends the proof of Theorem \ref{theo:main_u}.
\end{proof}

\vspace{5mm}
\noindent
{\bf \large Notation.}
\begin{itemize}
\item We always assume  by default that the summation index avoids the terms with vanishing denominators, and that the summation indices do not vanish. To illustrate this notation, consider the example
\[
\sum_{n_1+n_2+n_3=0} \frac{f(n_1,n_2,n_3)}{n_1(n_1-n_2)}: = \sum_{\stackrel{n_1+n_2+n_3=0}{n_1\neq 0, n_2\neq 0, n_3\neq 0, n_1\neq n_2}} \frac{f(n_1,n_2,n_3)}{n_1(n_1-n_2)}.
\]
\item The expressions under the square roots are always taken over the absolute values, {\em i.e.}
$\sqrt{f}:=\sqrt{|f|}$.
\item $\D_q F$ is the sequence $\frac{\D F}{\D q(-n)}$.
\item
We use $\les$ sign to avoid using unimportant constants: \\
$A\les B$ means there is an absolute constant $K$ such that $A\leq K B$. In some cases the constant will depend on parameters such as $s$.\\
$A\les B(\eta-)$ means that for any $\gamma>0$, $A\leq C_\gamma B(\eta-\gamma)$. \\
$A\les B(\eta+)$ is defined similarly.
\item $\la n \ra = \sqrt{1+n^2}$
\item We denote the $k$th derivative of $H$ over the flow of $F$ by $g_F^k H$, which is defined iteratively as follows
$$
g_F^0 H := H,\quad g^{}_F H = g_F^1 H =\{H,F\}, \quad g_F^k H := \{g_F^{k-1} H, F\},\,k=2,3,4,...
$$
\end{itemize}

\section{Canonical transformations}
The goal of this section is to transform the Hamiltonian to a more convenient (so called normal)
form where the most essential (resonant) terms are left at the low order. The non-resonant terms
will be absorbed into appropriate canonical transformations. Resonant terms are those that are
constant over the linear Hamiltonian flow, generated by  $\Lambda_2$, see the formal definition below.
In this section, the separation into the higher and lower order terms is formal as we will not invoke any estimates, yet. The results of this section are not new and follow closely
the standard normal form calculations, see {\em e.g.} \cite{ KappelerPoschel,Kuk}.

Consider the change of variables $u=u(q)$, generated by
the time-1 flow of a purely imaginary Hamiltonian $F$.
Namely, solve
\be \label{eq:Fflow}
\frac{dw(n)}{d\tau} = \sigma(n) \frac{\D F}{\D w(-n)},\,\,\,\,\,n\neq 0
\ee
with initial conditions
\[
w|_{\tau=0}=q,
\]
 thus producing a symplectic
transformation  $u = u(q):=\Phi_F^1(q)=w|_{\tau=1}$. On the other hand, we can write $q=\Phi_F^{-1}(u)$.
Let $\Phi_F^\tau$ be the time $\tau$ map of the flow of $F$. Using
Taylor expansion in $\tau$, evaluated at $\tau = 1$, we have
\begin{align}\label{eq:taylor}
H\circ \Phi_F^1(q)& =H(q) + \{H,F\}(q)+ \ldots
+ \frac{1}{k!}\{\ldots\{\{H,\underbrace{F\},F\},\ldots,F}_{k }\}(q)  \\
& + \int_0^1 \frac{(1-\tau)^k}{k!} \{\ldots\{\{H,\underbrace{F\},F\},\ldots,F}_{k+1}\}\circ \Phi_F^\tau(q) \,d\tau,   \nonumber
\end{align}
where the Poisson bracket is defined as the derivative of one Hamiltonian function, over
the flow of the other one
\be
\{A,B \} =  \sum_{n\neq 0}  \sigma(n) \frac{\D A}{\D q(n)} \frac{\D B}{\D q(-n)} .
\ee
Using the notation $g^k_F(H)$, see above, we can rewrite \eqref{eq:taylor} as
\begin{align}\label{eq:taylorg}
&H\circ \Phi_F^1 = \\
 &  H+g^{}_F H +\frac{1}{2}g_F^2 H +\ldots+\frac{1}{k!}g_F^k H  +\int_0^1 \frac{(1-\tau)^k}{k!} \big(g_F^{k+1} H \big)\circ \Phi_F^\tau \,d\tau. \nn
\end{align}

We choose the first transformation as a time-1 map of the Hamiltonian flow
of purely imaginary Hamiltonian function
\[
F_1 = \sum_{n_1+n_2+n_3=0} \cF_1 (n_1,n_2,n_3) u(n_1)u(n_2)u(n_3).
\]
With this choice of symplectic structure, all Hamiltonian functions must be purely imaginary.
In particular, this Hamiltonian function is purely imaginary provided
\[
\cF_1(-n_1,-n_2,-n_3)= \overline{\cF_1(n_1,n_2,n_3)}.
\]

Using \eqref{eq:taylorg} with $k=2$ we have
\begin{align}
H\circ \Phi_{F_1}^1 =H+g^{}_{F_1} H +\frac{1}{2}g_{F_1}^2 H +\frac{1}{2} \int_0^1  (1-\tau)^2 \big(g_{F_1}^{3} H \big)\circ \Phi_{F_1}^\tau \,d\tau. \label{eq:taylorF1}
\end{align}

\begin{defin}
The monomial $M_{n_1 n_2 ... n_k} = q(n_1)q(n_2)...q(n_k)$ is called resonant if it commutes with
the linear flow, i.e.
\be
\{\Lambda_2,M  \} = 0.
\ee
Otherwise, the monomial is called non-resonant. 
The sum of monomials is called resonant (non-resonant) if all monomials
are resonant (non-resonant). We will write $\{\,\, \} =\{ \,\, \}^{\rm r} + \{\,\,  \}^{\rm nr}$, 
where $\{ \,\, \}^{\rm r}$ represents  resonant terms and $\{ \,\, \}^{\rm nr}$ represents nonresonant 
terms.
\end{defin}

Rewriting the Hamiltonian with $H= \Lambda_2 + H_3$, we have
\begin{align} \nn
&H\circ \Phi^1_{F_1} =\\ \nn
& \Lambda_2 + H_3 + \{\Lambda_2, F_1 \} +  \{H_3,F_1 \} +
\frac{1}{2}\{\{\Lambda_2 ,F_1 \},F_1 \} + \frac{1}{2 }\{\{H_3,F_1 \},F_1 \} + R_1,
\end{align}
where
\be \label{eq:R1def}
R_1=\frac{1}{2} \int_0^1  (1-\tau)^2 \big(g_{F_1}^{3} H \big)\circ \Phi_{F_1}^\tau \,d\tau.
\ee
We choose $F_1$ so that to eliminate cubic non-resonant terms (in our case all
cubic terms are non-resonant)
\be
\{\Lambda_2, F_1\} +  H_3 = 0.
\label{eq:hom1}
\ee
Then we   have
\be
H\circ \Phi^1_{F_1}= \Lambda_2 + \frac{1}{2} \{H_3,F_1 \} + \frac{1}{2 }\{\{H_3,F_1 \},F_1 \} + R_1.
\ee
It turns out that another transformation $F_2$ that  removes non-resonant terms
in  $  \{H_3,F_1\}$ is required. For this purpose, we choose $F_2$ so that
\be \label{eq:F2def}
\{ \Lambda_2,F_2 \} + \frac{1}{2} \{H_3,F_1\} =\frac{1}{2} \{H_3,F_1\}^{\rm r} = \frac{3}{2} i \sum_{n\neq 0} |q(n)|^4,
\ee
see \eqref{eq:4res_term}.
Applying \eqref{eq:taylor} with $k=1$, the new Hamiltonian takes the form
\begin{align*}
H\circ \Phi^1_{F_1} \circ \Phi^1_{F_2}  = &\Lambda_2 + \frac{1}{2} \{H_3,F_1 \} + \frac{1}{2 }\{\{H_3,F_1 \},F_1 \} + R_1 \\
&+ \{ \Lambda_2,F_2 \} + \frac{1}{2} \{\{H_3,F_1\},F_2 \} + \frac{1}{2 }\{ \{\{H_3,F_1 \},F_1 \},F_2\} \\ & + \{R_1 ,F_2 \}
 + \int_0^1 (1-\tau) \{ \{H\circ \Phi^1_{F_1},F_2  \}, F_2 \} \circ \Phi^{\tau}_{F_2} d \tau.
\end{align*}
Using \eqref{eq:F2def},
we rewrite
\ba
H\circ \Phi^1_{F_1}\circ \Phi^1_{F_2} = \Lambda_2 + R,
\label{eq:trans_ham}
\ea
where
\begin{align}
R  = & \,\, i\, \frac{3}{2} \sum_{n\neq 0} |q(n)|^4  +  \frac{1}{2 }g_{F_1}^2 H_3 +  \frac{1}{2} g^{}_{F_2}g^{}_{F_1}H_3 + \frac{1}{2 } g_{F_2}^2 g^{}_{F_1}H_3   \nn \\
&+  R_1 + g^{}_{F_2} R_1 + \int_0^1 (1-\tau) g_{F_2}^2(H\circ \Phi^1_{F_1})  \circ \Phi^{\tau}_{F_2} d \tau, \label{eq:Rfinal}
\end{align}
and $R_1$ is given by \eqref{eq:R1def}.

\subsection{Calculation of $F_1$.}

Straightforward calculations give
\be
\{\Lambda_2, F_1 \} = i \sum_{n_1+n_2+n_3 =0} (n_1^3+n_2^3+n_3^3)
\cF_1(n_1,n_2,n_3) q(n_1) q(n_2) q(n_3).
\ee
Note that under the restriction $n_1+n_2+n_3 = 0$, the sum of cubes can be factored out
\[
n_1^3+n_2^3+n_3^3 = 3n_1n_2n_3.
\]
Thus, from \eqref{eq:hom1} we have
\be
\cF_1(n_1,n_2,n_3) = -\frac{\sigma(n_1n_2n_3)}{3\sqrt{n_1 n_2 n_3}},
\ee
whenever $n_1 n_2 n_3 \neq 0$. Otherwise $\cF(n_1,n_2,n_3)=0$.

\subsection{Calculation of $F_2$.}
We need to solve
\be
\{ \Lambda_2,F_2 \} + \frac{1}{2} \{H_3,F_1\}^{\rm nr} = 0,
\ee
but first we need to distinguish the resonant and nonresonant terms of $\{H_3,F_1\}$:
\[
\{H_3,F_1\}=\{H_3,F_1\}^{\rm nr}+\{H_3,F_1\}^{\rm r}.
\]
Recall,
\[
F_1 = -\sum_{n_1+n_2+n_3=0} \frac{\sigma(n_1n_2n_3)}{3\sqrt{n_1 n_2 n_3}} q(n_1) q(n_2) q(n_3)
\]
and
\[
H_3 (q)= i\sum_{n_1+n_2+n_3=0} \sqrt{n_1 n_2 n_3} \,\, q(n_1)q(n_2)q(n_3).
\]
Now, we compute
\begin{align*}
&\{H_3,F_1 \}  = \sum_{n\neq 0} \sigma(n)\frac{\D H_3}{\D q(n)} \frac{\D F_1}{\D q(-n)} \\
&=-i \sum_{n\neq 0} \sigma(n) \,\, 3\sum_{n_1+n_2+n=0}\sqrt{n_1 n_2 n} \, q(n_1)q(n_2) \,\,  \sum_{k_1+k_2-n=0}\frac{\sigma(k_1 k_2 n)}{ \sqrt{k_1 k_2 n}}\, q(k_1)q(k_2) \\
&=-3 i \sum_{\stackrel{n_1+n_2+n_3+n_4=0 }{n_1+n_2\neq 0}} \sqrt{\frac{n_1 n_2}{n_3 n_4}} \sigma(n_3 n_4) q(n_1) q(n_2) q(n_3) q(n_4).
\end{align*}
The resonant terms are the ones satisfying $n_1^3+n_2^3+n_3^3+n_4^3=0$.
Since we can rewrite (under the restriction $n_1+n_2+n_3+n_4=0 $)
$$
n_1^3+n_2^3+n_3^3+n_4^3=3(n_1+n_2)(n_1+n_3)(n_2+n_3),
$$
and $n_1+n_2\neq 0$, the resonant terms are the ones with $n_1+n_3=0$ or $n_2+n_3=0$.
Therefore, the nonresonant terms are
$$
\{H_3,F_1\}^{\rm nr}  =  -3i \sum_{\stackrel{n_1+n_3 \neq 0,n_2+n_3\neq 0, n_1+n_2\neq 0}{n_1+n_2+n_3+n_4=0}} \sqrt{\frac{n_1 n_2}{n_3 n_4}} \, \sigma(n_3 n_4) \, q(n_1) q(n_2) q(n_3) q(n_4).
$$
On the other hand, the resonant terms are
\begin{align}
\{H_3,F_1\}^{\rm r} & =  -3 i \sum_{\stackrel{n_1+n_3=0, n_1+n_2\neq 0}{n_1+n_2+n_3+n_4=0}}(\cdots) -3 i \sum_{\stackrel{n_2+n_3=0, n_1+n_3 \neq 0, n_1+n_2\neq 0}{n_1+n_2+n_3+n_4=0}}(\cdots)\nn\\
&= - 3i \sum_{n_1+n_2\neq 0} \sigma(n_1 n_2) |q(n_1)|^2 |q(n_2)|^2 - 3i \sum_{n_1\pm n_2\neq 0} \sigma(n_1 n_2) |q(n_1)|^2 |q(n_2)|^2.
\label{eq:res_terms}
\end{align}
The resonant terms can be simplified as follows
\[
(\ref{eq:res_terms})=-6 i \sum_{n_1 \pm n_2 \neq 0} \sigma(n_1 n_2) |q(n_1)|^2 |q(n_2)|^2 +3i \sum_{n\neq 0} |q(n)|^4 = 3i \sum_{n\neq 0} |q(n)|^4,
\]
since the first sum is equal to zero due to the cancellations
\[
\sigma(n_1 n_2)+ \sigma(-n_1 n_2) = 0
\]
and $|q(m)|=|q(-m)|$.

\noindent
Therefore, we have
\begin{align}
\{H_3,F_1\}^{\rm r} =  3i \sum_{n\neq 0} |q(n)|^4.  
\label{eq:4res_term}
\end{align}
Next, we solve
\[
\{ \Lambda_2,F_2 \} + \frac{1}{2} \{H_3,F_1\}^{\rm nr} = 0.
\]
By straightforward calculations,  taking $F_2$ of the form
\[
F_2 = \sum_{n_1+n_2+n_3+n_4=0} \cF_2 (n_1,n_2,n_3,n_4) q(n_1)q(n_2)q(n_3)q(n_4),
\]
we obtain
\begin{align*}
 &\{ \Lambda_2,F_2 \} = \\& = i\sum_{n_1+n_2+n_3+n_4=0} (n_1^3+n_2^3+n_3^3+n_4^3)\cF_2(n_1,n_2,n_3,n_4)q(n_1)q(n_2)q(n_3)q(n_4).
\end{align*}
Therefore, $\cF_2$ must satisfy the equality
\[
(n_1^3+n_2^3+n_3^3+n_4^3)\cF_2(n_1,n_2,n_3,n_4) - \frac{3}{2}\sqrt{\frac{n_1 n_2}{n_3 n_4}} \sigma(n_3 n_4) = 0,
\]
on the ''non-resonant set"
\ba \nn
NR_4 = \{ (n_1,n_2,n_3,n_4)\in \Z^4,n_1+n_2+n_3+n_4=0, \,\, n_1^3+n_2^3+n_3^3+n_4^3 \neq 0 \}.
\ea
Thus,
\be
\cF_2(n_1,n_2,n_3,n_4) = \frac{3}{2}\sqrt{\frac{n_1 n_2}{n_3 n_4}} \sigma(n_3 n_4) \frac{1}{n_1^3+n_2^3+n_3^3+n_4^3},
\ee
if $(n_1,n_2,n_3,n_4) \in NR_4$, and $\cF_2 = 0$ otherwise.

\section{Proof of Proposition~\ref{prop:transf}}
We start with a-priori estimates for the derivatives of $F_1$ and $F_2$. We need these estimates also in the subsequent sections.
Define the sequence
$$f_1(q_1,q_2)(n):=-\sum_{n_1+n_2+n=0}\frac{\sigma(n_1n_2n)}{\sqrt{n_1n_2n}}q_1(n_1)q_2(n_2)$$
so that
$$\partial_q F_1(n) = \frac{\partial F_1}{\partial q(-n)}=f_1(q,q)(n).$$
\begin{lem}\label{lem:f1} The following a-priori estimates hold,
\begin{align*}
\|f_1(q_1,q_2)\|_{\ell^2_{0-}}&\lesssim \|q_1\|_{\ell^2_{-1/2}}\|q_2\|_{\ell^2_{-1/2}},\\
 \|f_1(q_1,q_2)\|_{\ell^2_{\frac{1}{2}-}}&\lesssim \|q_1\|_{\ell^2_{-1/2}}\|q_2\|_{\ell^2},\\
 \|f_1(q_1,q_2)\|_{\ell^2_{1-}}&\lesssim \|q_1\|_{\ell^2}\|q_2\|_{\ell^2},\\
 \|f_1(q_1,q_2)\|_{\ell^2_{\frac{3}{2}}}&\lesssim \|q_1\|_{\ell^2_{1/2}}\|q_2\|_{\ell^2_{0+}}+\|q_1\|_{\ell^2_{0+}}\|q_2\|_{\ell^2_{1/2}}.
\end{align*}
\end{lem}

\vspace{5mm}

Now, define the sequence
\begin{align*}
&f_2(q_1,q_2,q_3)(n)= \partial_q F_2(n) = \frac{\partial F_2}{\partial q(-n)}=\\
&3\sum_{n+n_1+n_2+n_3=0} \frac{\sqrt{\frac{nn_1}{n_2n_3}}\sigma(n_2n_3)+\sqrt{\frac{n_1n_2}{n_3n}}\sigma(n_3n)}{n^3+n_1^3+n_2^3+n_3^3}
q_1(n_1)q_2(n_2)q_3(n_3).
\end{align*}
%Note that with this notation, we have
%$$\partial_q F_2(n) = \frac{\partial F_2}{\partial q(-n)}=f_2(q,q,q)(n).$$
\begin{lem}\label{lem:f2} For any permutation $(j_1,j_2,j_3)$ of $(1,2,3)$, and for any \mbox{$s\in [-1,1]$,}
we have
$$
\|f_2(q_1,q_2,q_3)\|_{\ell^2_s}\lesssim \|q_{j_1}\|_{\ell^2_s} \|q_{j_2}\|_{\ell^2_{0+}}\|q_{j_3}\|_{\ell^2_{0+}}.
$$
Moreover,
$$
\|f_2(q_1,q_2,q_3)\|_{\ell^2_{3/2}}\lesssim \sum \|q_{j_1}\|_{\ell^2_{3/2}} \|q_{j_2}\|_{\ell^2_{0+}}\|q_{j_3}\|_{\ell^2_{0+}},
$$
where the sum is taken over all permutations $(j_1,j_2,j_3)$ of $(1,2,3)$.
\end{lem}
\vspace{5mm}

Now, we prove Proposition \ref{prop:transf} using Lemma~\ref{lem:f1} and Lemma~\ref{lem:f2}.
The proof of these lemmas will be given in the next section.

\begin{proof}[ Proof  of Proposition \ref{prop:transf}]
It suffices to prove that $\Phi_F^1$ is near identity for each $F=F_1$  and $F=F_2$ in the sense
 $$
q\in X^{\rho}_{\eps} \implies  \|\Phi_{F}^1(q)-q\|_{\ell^2_s}\lesssim \eps^{1-s-},\,\,\,s\in[0,3/2].
$$
This is because $\|\Phi_{F}^1(q)-q\|_{\ell^2_s}\lesssim \eps^{1-s-}$ implies that $\Phi_{F}^1(q)\in X^{2\rho}_{\eps}$ and because
if $\Phi_{F_1}^1$ and $\Phi_{F_2}^1$ are near identity, then their composition,
$\Phi_{F_1}^1 \circ \Phi_{F_2}^1$,
is also near identity.

Note that in light of equation \eqref{eq:Fflow} we have
\begin{align*}
\|\Phi_{F_1}^1(q)-q\|_{\ell^2_s}&=\Big\|\int_0^1 \sigma(n)\frac{\, \D F_1}{\D w(-n)} d\tau \Big\|_{\ell^2_s} \leq \Big\| \frac{\D F_1}{\D w(-n)} \Big\|_{\ell^2_s} =\|f_1(w,w)\|_{\ell^2_s}.
\end{align*}
Applying Lemma~\ref{lem:f1} with $q_1=q_2=w=\Phi_{F_1}^\tau(q)\in X^{2\rho}_{\eps}$,
we obtain
\ba
\|f_1(w,w)\|_{l^2} &\les& \eps \nn \\
\|f_1(w,w)\|_{l^2_{\frac32}} &\les& \eps^{\frac12-}, \nn
\ea
which implies that $\Phi_{F_1}^1$ is near identity.

Similarly, applying Lemma~\ref{lem:f2}, with $q_1=q_2=q_3 = w \in X^{\rho}_{\eps}$ we have
\ba
\|f_2(w,w,w)\|_{l^2} \les \eps^{\frac32-} \nn \\
\|f_2(w,w,w)\|_{l^2_{\frac32}} \les \eps^{0-}, \nn
\ea
which implies that $\Phi_{F_2}^1$ is near identity.

Since \eqref{eq:Fflow} is time reversible,   $\Phi_{F_1}^{-1}$ and $\Phi_{F_2}^{-1}$ are also near identity,
which implies that $q(u)\in X_\eps^{2\rho}$ if $u\in X_\eps^{\rho}$.
\end{proof}

\subsection{Proof of Lemma~\ref{lem:f1} and Lemma~\ref{lem:f2}}
We use the following lemma repeatedly in the proof of Lemma~\ref{lem:f1}, Lemma~\ref{lem:f2}, and in the subsequent sections. The proof is left to the reader.
\begin{lem}\label{lem:calc}
\vskip 0.1 in
a) For any $s, r \in \Bbb R$, $1 \leq p, q \leq \infty$ and $\frac{1}{p}-\frac{1}{q}>s-r \geq 0$, we have the embedding,
$$\|u\|_{l_{r}^{p}} \leq C\|u\|_{l_s^{q}}.$$
b) Let $1 \leq p, q, r \leq \infty$ and $\frac{1}{p}+\frac{1}{q}=1+\frac{1}{r}$, we have Young's inequality,
$$\|u*v\|_{l^{r}} \leq \|u\|_{l^{p}}\|v\|_{l^{q}}.$$
\end{lem}
For the convenience of the reader we record the definition of the discrete convolution
$$u*v(n)=\sum_{m }u(m)v(n-m)=\sum_{m  }u(n-m)v(m).$$
By a slight abuse of notation we  also denote by $u*v$ all the sums of the form $\sum_{m  }u(m)v(\pm n-m)$. Young's inequality holds true for all these
convolution products of functions.

\begin{proof}[Proof of Lemma~\ref{lem:f1}]
We begin with the second estimate. Note that 
$$
|f_1|\lesssim \frac{1}{\sqrt{n}}\big(\frac{|q_1|}{\sqrt{\cdot}} * \frac{|q_2|}{\sqrt{\cdot}}\big)(n).
$$
Therefore, for $s<1/2$,
\begin{align*}
 \|f_1\|_{\ell^2_s}&\lesssim \big\|\frac{|q_1|}{\sqrt{\cdot}} * \frac{|q_2|}{\sqrt{\cdot}}\big\|_{\ell^2_{s-1/2}} \lesssim \big\|\frac{|q_1|}{\sqrt{\cdot}} * \frac{|q_2|}{\sqrt{\cdot}}\big\|_{\ell^{2+}} \\
&\lesssim \big\|\frac{|q_1|}{\sqrt{\cdot}} \big\|_{\ell^2} \big\|\frac{|q_2|}{\sqrt{\cdot}}\big\|_{\ell^{1+}} \lesssim \big\|  q_1   \big\|_{\ell^2_{-1/2}} \big\| q_2 \big\|_{\ell^{2}}.
\end{align*}
The second inequality follows from the first part of Lemma~\ref{lem:calc}, and  the third inequality follows from Young's inequality. Finally the last inequality is another application of the
first part of the aforementioned Lemma. \\

To prove the first estimate, note that for $s<0$,
\begin{align*}
 \|f_1\|_{\ell^2_s}&\lesssim \big\|\frac{|q_1|}{\sqrt{\cdot}} * \frac{|q_2|}{\sqrt{\cdot}}\big\|_{\ell^2_{s-1/2}} \lesssim \big\|\frac{|q_1|}{\sqrt{\cdot}} * \frac{|q_2|}{\sqrt{\cdot}}\big\|_{\ell^{\infty}}\\
&\lesssim \big\|\frac{|q_1|}{\sqrt{\cdot}} \big\|_{\ell^2} \big\|\frac{|q_2|}{\sqrt{\cdot}}\big\|_{\ell^{2}} = \big\|  q_1   \big\|_{\ell^2_{-1/2}} \big\| q_2 \big\|_{\ell^{2}_{-1/2}}.
\end{align*}
Again the second inequality follows from the first part of Lemma \ref{lem:calc} and the third inequality follows from Young's inequality. \\

For the fourth estimate note that for $s>1/2$, using $|n|^{s-1/2}\lesssim |n_1|^{s-1/2}+|n_2|^{s-1/2}$, we have
$$|n|^s |f_1(n)|\lesssim  ( |\cdot|^{s-1}|q_1| ) * ( |\cdot|^{-1/2} |q_2|  )(n) + ( |\cdot|^{-1/2}|q_1| ) * ( |\cdot|^{s-1} |q_2|  )(n).
$$
Therefore, for $s>1/2$,
\begin{align*}
 \|f_1\|_{\ell^2_s}&\lesssim \big\|(|\cdot|^{s-1}|q_1| ) * ( |\cdot|^{-1/2} |q_2| )\big\|_{\ell^2} + \big\| (|\cdot|^{-1/2}|q_1| ) * ( |\cdot|^{s-1} |q_2|) \big\|_{\ell^2}.
\end{align*}
In particular, for $s=3/2$, we have
\begin{align*}
 \|f_1\|_{\ell^2_{3/2}}&\lesssim \big\|(|\cdot|^{1/2}|q_1| ) * ( |\cdot|^{-1/2} |q_2| )\big\|_{\ell^2} + \big\| (|\cdot|^{-1/2}|q_1| ) * ( |\cdot|^{1/2} |q_2|) \big\|_{\ell^2}\\
&\lesssim \big\| |\cdot|^{1/2} q_1 \big\|_{\ell^2} \big\| |\cdot|^{-1/2} q_2 \big\|_{\ell^1}+\big\| |\cdot|^{1/2} q_2 \big\|_{\ell^2} \big\| |\cdot|^{-1/2} q_1\big\|_{\ell^1} \\
&\lesssim \big\|  q_1 \big\|_{\ell^2_{1/2}} \big\| q_2 \big\|_{\ell^2_{0+}}+\big\|  q_2 \big\|_{\ell^2_{1/2}} \big\| q_1 \big\|_{\ell^2_{0+}}.
\end{align*}
For this estimate we first apply Young's inequality and then the first part of Lemma \ref{lem:calc}. \\

Finally to prove the third estimate, for fixed $\delta>0$
\begin{align*}
 \|f_1\|_{\ell^2_{1-\delta}}&\lesssim \big\|(|\cdot|^{-\delta}|q_1| ) * ( |\cdot|^{-1/2} |q_2| )\big\|_{\ell^2} + \big\| (|\cdot|^{-1/2}|q_1| ) * ( |\cdot|^{-\delta} |q_2|) \big\|_{\ell^2}\\
&\lesssim \big\| |\cdot|^{-\delta} q_1 \big\|_{\ell^{2-\eta}} \big\| |\cdot|^{-1/2} q_2 \big\|_{\ell^{1+\tilde\eta }}+\big\| |\cdot|^{-\delta} q_2 \big\|_{\ell^{2-\eta}} \big\| |\cdot|^{-1/2} q_1 \big\|_{\ell^{1+\tilde\eta }}\\
&\lesssim  \|q_1\|_{\ell^2}\|q_2\|_{\ell^2}.
\end{align*}
The derivation of this last string of inequalities follows as above using Lemma \ref{lem:calc}. We only note that the last step follows if we choose $\eta, \tilde\eta >0$ sufficiently small with $\eta$
depending on $\delta>0$.
\end{proof}

\begin{proof}[Proof of Lemma~\ref{lem:f2}]
First note that by duality and interpolation it suffices to prove the first statement for $s=1$.
To estimate $\|f_2\|_{\ell^2_1}$ we use duality as follows
$$ \|f_2\|_{\ell^2_1} = \sup_{\|h\|_{\ell^2_{-1}=1}} |\langle f_2,h \rangle|.$$
Note that for any permutation $(j_1,j_2,j_3)$ the form $\langle f_2,h \rangle$ on the right hand side can be estimated by
\begin{multline} \nn
\sum_{n_1+n_2+n_3+n_4=0} \frac{|n_1n_2|+|n_1n_3|+|n_1n_4|+|n_2n_3|+|n_2n_4|+|n_3n_4|}{\sqrt{n_1n_2n_3n_4}|n_1^3+n_2^3+n_3^3+n_4^3|} \times\\ \times
|q_{j_1}(n_1)q_{j_2}(n_2)q_{j_3}(n_3)h(n_4)|.
\end{multline}
 Since $q_{j_2}$ and $q_{j_3}$ appear symmetrically on the right side of the inequality, it suffices to estimate  the following sum
$$
 \sum_{n_1+n_2+n_3+n_4=0} \frac{|n_1n_2|+  |n_1n_4|+|n_2n_3|+|n_2n_4|  }{\sqrt{n_1n_2n_3n_4}|n_1^3+n_2^3+n_3^3+n_4^3|}
 |q_{j_1}(n_1)q_{j_2}(n_2)q_{j_3}(n_3)h(n_4)|.
$$
The estimate for these summands are very similar to each other, therefore we  consider only two of them. The remaining estimates just recycle the arguments below and will be omitted:
$$
\sum_{n_1+n_2+n_3+n_4=0} \frac{ |n_2n_3|+|n_2n_4| }{\sqrt{n_1n_2n_3n_4}|n_1^3+n_2^3+n_3^3+n_4^3|}
|q_{j_1}(n_1)q_{j_2}(n_2)q_{j_3}(n_3)h(n_4)|
$$
\begin{multline} \label{eq:f21}
 =\sum_{n_1+n_2+n_3+n_4=0} \frac{ |n_2|^{1/2-\delta}  |n_3|^{ 1/2-\delta}  |n_4|^{1/2}  }{|n_1|^{3/2} \, |n_1^3+n_2^3+n_3^3+n_4^3|} \times \\ \times
\big| n_1q_{j_1}(n_1) |n_2|^\delta q_{j_2}(n_2)|n_3|^\delta q_{j_3}(n_3)\frac{h(n_4)}{n_4} \big|
\end{multline}
\begin{multline}
+\sum_{n_1+n_2+n_3+n_4=0} \frac{ |n_2|^{1/2-\delta}   |n_4|^{3/2} }{|n_1|^{3/2} |n_3|^{1/2+\delta} \, |n_1^3+n_2^3+n_3^3+n_4^3|} \times \\ \times
\big| n_1q_{j_1}(n_1) |n_2|^\delta q_{j_2}(n_2)|n_3|^\delta q_{j_3}(n_3)\frac{h(n_4)}{n_4} \big|. \label{eq:f22}
\end{multline}
After substituting $n_4=-n_1-n_2-n_3$, the multiplier in \eqref{eq:f22} takes the form
\begin{align*}
&\frac{ |n_2|^{1/2-\delta} \, |n_1+n_2+n_3|^{3/2} }{|n_1|^{3/2}  |n_3|^{1/2+\delta}\, |n_1 +n_2\| n_1+n_3| |n_2+n_3|}\\ &\lesssim \frac{|n_2|^{1/2-\delta}  }{   |n_3|^{1/2+\delta} \, |n_1 +n_2\| n_1+n_3| |n_2+n_3|}+\frac{|n_2|^{1/2-\delta} |n_2+n_3|^{1/2} }{|n_1|^{3/2}  |n_3|^{1/2+\delta} \, |n_1 +n_2\| n_1+n_3| }.
\end{align*}
\vskip 0.1 in
Note that  $|n_2|\lesssim |n_1 +n_2\| n_1+n_3| |n_2+n_3| $ and $|n_2|\lesssim |n_1\|n_1+n_2|$  under the condition $|n_1n_2n_3| |n_1 +n_2\| n_1+n_3| |n_2+n_3|\neq 0$. Using this we further bound the multiplier by
\begin{align}
 & \frac{  1 }{ |n_2 n_3|^{1/2+\delta}}  +\frac{ |n_2+n_3|^{1/2} }{\sqrt{n_1}|n_2 n_3|^{1/2+\delta} \, | n_1+n_3| } \nonumber\\
&\lesssim \frac{  1 }{|n_2 n_3|^{1/2+\delta}}  +\frac{1}{\sqrt{n_1}   |n_3|^{1/2+\delta} \, |n_1+n_3| }+\frac{1}{\sqrt{n_1}   |n_2 |^{1/2+\delta} | n_3|^{ \delta}\, | n_1+n_3| } \nonumber
\end{align}
\begin{align}\lesssim \frac{  1 }{|n_2 n_3|^{1/2+\delta}}  +\frac{1}{|n_1 n_3|^{1/2+\delta/2}}+\frac{  1 }{|n_1 n_2|^{1/2+\delta}}. \label{eq:f23}
\end{align}
\vskip 0.1 in
The last inequality follows from  $|n_1|\lesssim |n_3\|n_1+n_3|$.
The contribution of the first summand in \eqref{eq:f23} to \eqref{eq:f22} is
\begin{align*}
& \sum_{n_1+n_2+n_3+n_4=0}  \frac{  1 }{|n_2 n_3|^{1/2+\delta}}
\big| n_1q_{j_1}(n_1) |n_2|^\delta q_{j_2}(n_2)|n_3|^\delta q_{j_3}(n_3)\frac{h(n_4)}{n_4} \big| \\
& =\sum_{n_1,n_2,n_3}  \frac{  1 }{|n_2 n_3|^{1/2+\delta}}
\big| n_1q_{j_1}(n_1) |n_2|^\delta q_{j_2}(n_2)|n_3|^\delta q_{j_3}(n_3)\frac{h(-n_1-n_2-n_3)}{ n_1+n_2+n_3} \big|\\
&\lesssim \|q_{j_1}\|_{\ell^2_1}  \|h\|_{\ell^2_{-1}}\sum_{n_2,n_3}  \frac{  1 }{|n_2 n_3|^{1/2+\delta}}
\big|  |n_2|^\delta q_{j_2}(n_2)|n_3|^\delta q_{j_3}(n_3) \big|\\
&\lesssim \|q_{j_1}\|_{\ell^2_1}\|q_{j_2}\|_{\ell^2_\delta}\|q_{j_3}\|_{\ell^2_\delta}\|h\|_{\ell^2_{-1}}.
\end{align*}
The first inequality follows from Cauchy-Schwarz in $n_1$ sum and the second follows from Cauchy-Schwarz in
$n_2, n_3 $ sums since $\frac{1}{|n_2 n_3|^{1/2+\delta}}$ is square summable.  The contribution of the other two
summands in \eqref{eq:f23} to \eqref{eq:f22} can be estimated similarly.

Now we consider \eqref{eq:f21}. After
substituting $n_4=-n_1-n_2-n_3$, the multiplier takes the form
\begin{align*}
&\frac{ |n_2n_3|^{1/2-\delta} \, |n_1+n_2+n_3|^{1/2} }{|n_1|^{3/2}  \, |n_1 +n_2\| n_1+n_3| |n_2+n_3|}\\ &\lesssim \frac{|n_2n_3|^{1/2-\delta}  }{   |n_1|  \, |n_1 +n_2\| n_1+n_3| |n_2+n_3|}+\frac{|n_2n_3|^{1/2-\delta}  }{|n_1|^{3/2}   \, |n_1 +n_2\| n_1+n_3| |n_2+n_3|^{1/2}}\\
&\lesssim \frac{1 }{ |n_1|^{2\delta}|n_1 +n_2|^{1/2+\delta}| n_1+n_3|^{1/2+\delta} |n_2+n_3|}+\frac{1  }{|n_2 |^{1/2+\delta}|n_3|^\delta |n_2+n_3|^{1/2}}\\
&\lesssim \frac{1 }{  | n_2-n_3|^{1/2+\delta}   |n_2+n_3|}+\frac{1  }{|n_2 |^{1/2+\delta}|n_3|^\delta |n_2+n_3|^{1/2}}.
\end{align*}

We estimate the two terms separately. By summing first in $n_1$ and then using Cauchy-Schwarz inequality in $n_2, n_3$, to estimate
\begin{multline}\nn
\sum_{n_1,n_2,n_3}  \frac{  1 }{|n_2-n_3|^{1/2+\delta}|n_2+n_3|} \times \\ \times
\big| n_1q_{j_1}(n_1) |n_2|^\delta q_{j_2}(n_2)|n_3|^\delta q_{j_3}(n_3)\frac{h(-n_1-n_2-n_3)}{ n_1+n_2+n_3} \big|
\end{multline}
it is enough to bound
$$\sum_{n_2,n_3}\frac{1}{|n_2-n_3|^{1+2\delta} |n_2+n_3|^2}.$$
But
$$\sum_{n_2,n_3}\frac{1}{|n_2-n_3|^{1+2\delta} |n_2+n_3|^2} = \sum_{n_2,m}\frac{1}{|m|^{1+2\delta} |2n_2+m|^2} <\infty.$$
\\
\\
To estimate the second term by the above arguments it is enough to estimate
$$\sum_{n_2,n_3}\frac{1}{|n_2|^{1+2\delta}}\frac{1}{|n_3|^{2\delta}|n_2+n_3|}.$$
But
$$\sum_{n_2,n_3}\frac{1}{|n_2|^{1+2\delta}}\frac{1}{|n_3|^{2\delta}|n_2+n_3|}=\sum_{n_2}\frac{1}{|n_2|^{1+2\delta}}\left ( \sum_{n_3}\frac{1}{|n_3|^{2\delta}|n_2+n_3|}\right )$$
and H\"older's inequality implies that
$$\|n_{3}^{-2\delta}(n_2+n_3)^{-1}\|_{l^{1}(n_3)} \lesssim \|n_{3}^{-2\delta}\|_{l_{n_3}^{\infty-}}\|\frac{1}{n_2+n_3}\|_{l^{1+}(n_3)} <\infty,$$
while the rest is summable in $n_2$. This finishes the proof of the first assertion of the lemma.
\\
\\
To prove the second assertion we use duality in a similar way. Since we have a sum over all possible permutations in the right hand side of the inequality it suffices to consider the  following sum
$$\sum_{n_1+n_2+n_3+n_4=0} \frac{|n_1n_3|+|n_1n_4| }{\sqrt{n_1n_2n_3n_4}|n_1^3+n_2^3+n_3^3+n_4^3|}
 |q_{j_1}(n_1)q_{j_2}(n_2)q_{j_3}(n_3)h(n_4)|$$
\begin{multline}\nn
=\sum_{n_1+n_2+n_3+n_4=0} \frac{  |n_3|^{ 1/2-\delta}  |n_4|  }{|n_1| |n_2|^{1/2+\delta}  \, |n_1^3+n_2^3+n_3^3+n_4^3|} \times \\ \times
\big| |n_1|^{3/2}q_{j_1}(n_1) |n_2|^\delta q_{j_2}(n_2)|n_3|^\delta q_{j_3}(n_3)\frac{h(n_4)}{|n_4|^{3/2}} \big|
\end{multline}
\begin{multline}\nn +\sum_{n_1+n_2+n_3+n_4=0} \frac{   |n_4|^{ 2} }{|n_1|  |n_2n_3|^{1/2+\delta}  \, |n_1^3+n_2^3+n_3^3+n_4^3|} \times \\ \times
\big| |n_1|^{3/2}q_{j_1}(n_1) |n_2|^\delta q_{j_2}(n_2)|n_3|^\delta q_{j_3}(n_3)\frac{h(n_4)}{|n_4|^{3/2}} \big|.
\end{multline}

As above the proof follows from the following estimate for the multipliers:
\begin{multline} \nn
\frac{  |n_3|^{ 1/2-\delta}  |n_4|  }{|n_1| |n_2|^{1/2+\delta}  \, |n_1^3+n_2^3+n_3^3+n_4^3|}
+ \frac{   |n_4|^{ 2} }{|n_1|  |n_2n_3|^{1/2+\delta}  \, |n_1^3+n_2^3+n_3^3+n_4^3|}\\ \lesssim \frac{1}{|n_2n_3|^{1/2+\delta}}.
\end{multline}
To prove this inequality we first substitute  $n_4=-n_1-n_2-n_3$ to obtain
\begin{multline}
\frac{  |n_3|^{ 1/2-\delta}  |n_1+n_2+n_3|  }{|n_1| |n_2|^{1/2+\delta}  \, |n_1 +n_2\| n_1+n_3| |n_2+n_3|}
+ \\ + \frac{   |n_1+n_2+n_3|^{ 2} }{|n_1|  |n_2n_3|^{1/2+\delta}  \,|n_1 +n_2\| n_1+n_3| |n_2+n_3|}\nonumber
\end{multline}
\begin{align}
&\lesssim \frac{  |n_3|^{ 1/2-\delta}     }{  |n_2|^{1/2+\delta}  \, |n_1 +n_2\| n_1+n_3| |n_2+n_3|}
+\frac{  |n_3|^{ 1/2-\delta}    }{|n_1| |n_2|^{1/2+\delta}  \, |n_1 +n_2\| n_1+n_3| }\nonumber\\
&+ \frac{   |n_1|}{  |n_2n_3|^{1/2+\delta}  \,|n_1 +n_2\| n_1+n_3| |n_2+n_3|}+\frac{   | n_2+n_3|  }{|n_1|  |n_2n_3|^{1/2+\delta}  \,|n_1 +n_2\| n_1+n_3|  }.  \nn
\end{align}
Using the inequalities
\begin{align*}
|n_3|,|n_1|&\lesssim |n_1 +n_2\| n_1+n_3| |n_2+n_3| \\
|n_3| &\lesssim |n_1\|n_1+n_3|\\
|n_2+n_3|&\leq |n_2|+|n_3|\lesssim |n_1\|n_1+n_2|+|n_1\|n_1+n_3|
\end{align*}
we see that last line is bounded by $\frac{1}{|n_2n_3|^{1/2+\delta}}$. This finishes the proof of the lemma by using the methods of the first part.
\end{proof}

\section{Proof of Proposition~\ref{prop:error}}
Note that by \eqref{eq:Eqn} and \eqref{eq:Rfinal}, it suffices to prove that if $q\in X_\eps^\rho$, then,   for $s\in[0,1/2]$,
\begin{align} \label{eq:qn}
\|q(k)^3\|_{\ell^2_s}&\lesssim \eps^{1-s-}\\
\label{eq:gfabest} \Big\|\D_q \, g_{F_2}^b g_{F_1}^a H_3 \Big\|_{\ell^2_s}&\lesssim \eps^{1-s-},\,\,\,\,\,\text{ if } a\geq 1, a+b\geq 2,
\end{align}
and similarly for terms involving integrals.

The inequality \eqref{eq:qn} is obtained as follows:
$$
\|q^3\|_{\ell^2_s}\leq \|q\|_{\ell^\infty}^2 \|q \|_{\ell^2_s} \lesssim \eps \, \eps^{\frac12-s} =\eps^{\frac32-s}.
$$
The inequality \eqref{eq:gfabest} follows from Theorem~\ref{thm:H3ab} in the next section, and the estimates for the integral terms are discussed in  section~\ref{sec:rmn}.

\subsection{Estimates for the terms $g_{F_2}^b g_{F_1}^a H_3 $}
In this section, we estimate the derivative of the commutators $g_{F_2}^bg_{F_1}^a H_3$ for $a\geq 1$.
\begin{theo}\label{thm:H3ab}
Let $\eps>0$. Assume that $q\in X_\rho^\eps$.  Then
for $a\geq 1$, $b\geq 0$, and $s\in[0,1/2],$ we have
$$\Big\|\D_q \, g_{F_2}^bg_{F_1}^a H_3\Big\|_{\ell^2_s}\lesssim \eps^{\frac{a}{2}+b-s-}.$$
\end{theo}
With a slight abuse of notation, define
$$H_3(q_1,q_2,q_3)=i \sum_{n_1+n_2+n_3=0} \sqrt{n_1 n_2 n_3} \,\, q_1(n_1)q_2(n_2)q_3(n_3).$$
With this notation, $H_3(q)=H_3(q,q,q)$.
Note that $\{H_3,F_1\}$ can be written as
\begin{align}\label{H3F1}
 H_3( \partial_q F_1, q,q)+H_3( q,\partial_q F_1, q)+H_3(q,q, \partial_q F_1),
\end{align}
where $\partial_q F_1$ is the sequence $\frac{\partial F_1}{\partial q(-n)}=f_1(q,q)(n)$. By symmetry, we can write
 $$\{H_3,F_1\}=3 \, H_3( f_1(q,q), q,q)$$
By the same token, we can write $\{\{H_3,F_1\},F_2\}$ as a sum of the following terms
$$H_3(f_1(q,q),f_2(q,q,q),q),\,\,\text{and } H_3(f_1(f_2(q,q,q),q),q,q).$$
To generalize this to higher order commutators, we define $Q_{a,b}$ as follows. First $Q_{0,0}$ is $q$.
To obtain $Q_{a,b}$, start with $Q_{0,0}$ and iteratively, a times, replace one $q$ with $f_1(q,q)$, then again iteratively replace one $q$ with $f_2(q,q,q)$ $b$ times. Any sequence obtained in this manner is called $Q_{a,b}$.
For example $\{\{H_3,F_1\},F_2\}$ can be described as a sum of
$$H_3(Q_{1,0},Q_{0,1},Q_{0,0}),\,\,\text{and } H_3(Q_{1,1},Q_{0,0},Q_{0,0}).$$
In general, we can write $g_{F_2}^bg_{F_1}^a H_3$ as a sum of terms of the form
\begin{align}\label{H3ab}
H_3( Q_{a_1,b_1}, Q_{a_2,b_2}, Q_{a_3,b_3}), \,\,a_1+a_2+a_3=a, \,\,b_1+b_2+b_3=b, a_j, b_j \in \mathbb N.
\end{align}
To estimate $\|\partial_q g_{F_2}^bg_{F_1}^a H_3\|_{\ell^2_s}$, we use duality and estimate
$$\sup_{\|h\|_{\ell^2_{-s}}=1}|\langle \partial_q g_{F_2}^bg_{F_1}^a H_3 , h \rangle |.
$$
Note that $\langle \partial_q g_{F_2}^bg_{F_1}^a H_3 , h \rangle$ can be written as the commutator $\{g_{F_2}^bg_{F_1}^a H_3,G\}$, where $G=G(q)=\sum_n h(n)q(n)$. This is because $h=\partial_q G$.
In light of \eqref{H3F1} and \eqref{H3ab}, we can now write $\langle \partial_q g_{F_2}^bg_{F_1}^a H_3 , h \rangle$ as a sum of terms of the form
\be\label{eq:H3qabh}
H_3( Q_{a_1,b_1}^h, Q_{a_2,b_2}, Q_{a_3,b_3}), \,\,a_1+a_2+a_3=a, \,\,b_1+b_2+b_3=b, a_j, b_j \in \mathbb N,
\ee
where $Q_{a,b}^h$ is obtained from $Q_{a,b}$ by replacing one $q$ by $h$.
Therefore, to prove Theorem~\ref{thm:H3ab}, we need to estimate the sequences $Q_{a,b}$ and $Q_{a,b}^h$ in $\ell^2_s$ spaces, and estimate $H_3(q_1,q_2,q_3)$ for $q_j$ in $\ell^2_s$ spaces:
\begin{prop}\label{prop:H3} For any permutation $(j_1,j_2,j_3)$ of $(1,2,3)$, we have
\begin{align*}
|H_3(q_1,q_2,q_3)|&\lesssim \|q_{j_1}\|_{\ell^2} \| q_{j_2}\|_{\ell^2_{1+}} \|q_{j_3}\|_{\ell^2_{1+}} ,\\
|H_3(q_1,q_2,q_3)|&\lesssim \|q_{j_1}\|_{\ell^{2}_{1-}} \| q_{j_2}\|_{\ell^2_{1/2+}} \|q_{j_3}\|_{\ell^2_{1/2+}},\\
|H_3(q_1,q_2,q_3)|&\lesssim \|q_{j_1}\|_{\ell^2_{-1/2}} \big(\| q_{j_2}\|_{\ell^2_{3/2}} \|q_{j_3}\|_{\ell^2_{1+}}+\| q_{j_2}\|_{\ell^2_{1+}} \|q_{j_3}\|_{\ell^2_{3/2}}\big),\\
|H_3(q_1,q_2,q_3)|&\lesssim \|q_{j_1}\|_{\ell^{2}_{1/2-}} \| q_{j_2}\|_{\ell^2_{1 }} \|q_{j_3}\|_{\ell^2_{1/2+}}
\end{align*}
\end{prop}

\begin{prop} \label{prop:Qab}
Assume that $q$ satisfies the hypothesis of Theorem~\ref{thm:H3ab}. Then we have
\begin{align*}
\|Q_{a,b}\|_{\ell^2_s}&\lesssim \eps^{\frac32-s+\frac{a}{2}+b-},\,\,\,s\in[1,3/2],\,\,a\geq 1,\,\,b\geq 0\\
\|Q_{0,b}\|_{\ell^2_s}&\lesssim \eps^{\frac12-s+b-},\,\,\,s\in[0,3/2], \,\,b\geq 0.
\end{align*}
\end{prop}
\begin{prop}\label{prop:Qabh}
Assume that $q$ satisfies the hypothesis of Theorem~\ref{thm:H3ab}. Then for any $a\geq 1$, $b \geq 0$, we have
\begin{align*}
\|Q_{a,b}^h\|_{\ell^{2}_{1-}}&\lesssim \|h\|_{\ell^2}\,\eps^{\frac{a}{2}+b-}, \\
\|Q_{a,b}^h\|_{\ell^{2}_{1/2-}}&\lesssim \|h\|_{\ell^2_{-1/2}}\,\eps^{\frac{a}{2}+b-},
\end{align*}
and for any $b\geq 0$, we have
\begin{align*}
\|Q_{0,b}^h\|_{\ell^2}&\lesssim \|h\|_{\ell^2}\,\eps^{b-}, \\
\|Q_{0,b}^h\|_{\ell^{2}_{-1/2}}&\lesssim \|h\|_{\ell^2_{-1/2}}\,\eps^{b-}.
\end{align*}
\end{prop}
We now prove Theorem~\ref{thm:H3ab} using Propositions~\ref{prop:H3}, \ref{prop:Qab} and \ref{prop:Qabh}.
\begin{proof}[Proof of Theorem~\ref{thm:H3ab}]
By the discussion leading to \eqref{eq:H3qabh} we need to estimate $H_3( Q_{a_1,b_1}^h, Q_{a_2,b_2}, Q_{a_3,b_3})$ for $h\in \ell^2$ and $h\in \ell^2_{-1/2}$. First consider the case $h\in \ell^2$. We have the following subcases $a_1=0$, and $a_1\neq 0$.
In the former case, by Proposition~\ref{prop:H3} and Proposition~\ref{prop:Qabh}, we  have
\begin{align*}
|H_3(Q_{0,b_1}^h, Q_{a_2,b_2}, Q_{a_3,b_3})|&\lesssim   \|Q_{0,b_1}^h\|_{\ell^2} \| Q_{a_2,b_2}\|_{\ell^2_{1+}} \|Q_{a_3,b_3}\|_{\ell^2_{1+}}\\
&\lesssim \|h\|_{\ell^2} \, \eps^{b_1-} \| Q_{a_2,b_2}\|_{\ell^2_{1+}} \|Q_{a_3,b_3}\|_{\ell^2_{1+}}.
\end{align*}
Now, by Proposition~\ref{prop:Qab}, it is easy to see that the worst case is when $a_2=a$, $a_3=0$, in which case we obtain
\begin{align*}
|H_3(Q_{0,b_1}^h, Q_{a_2,b_2}, Q_{a_3,b_3})|&\lesssim  \|h\|_{\ell^2} \, \eps^{b_1-}  \eps^{\frac12+\frac{a}{2}+b_2-}\eps^{-\frac12+b_3-}\\
&\lesssim \|h\|_{\ell^2} \,\eps^{\frac{a}{2}+b-}.
\end{align*}
If $a_1\neq 0$, the worst case is when $a_1=a, a_2=a_3=0$. Using the Propositions above we have
\begin{align*}
|H_3(Q_{a,b_1}^h, Q_{0,b_2}, Q_{0,b_3})|&\lesssim   \|Q_{a,b_1}^h\|_{\ell^{2}_{1-}} \| Q_{0,b_2}\|_{\ell^2_{1/2+}} \|Q_{0,b_3}\|_{\ell^2_{1/2+}}\\
&\lesssim \|h\|_{\ell^2} \,\eps^{\frac{a}{2}+b_1-} \eps^{b_2-}\eps^{b_3-}\\
&\lesssim \|h\|_{\ell^2} \,\eps^{\frac{a}{2}+b-}.
\end{align*}
It remains to consider the case $h\in \ell^2_{-1/2}$. As before we have the subcases $a_1=0$, $a_1\neq 0$.
If $a_1=0$, the worst case is when $a_2=a, a_3=0$. We estimate
\begin{align*}
& |H_3(Q_{0,b_1}^h, Q_{a ,b_2}, Q_{0,b_3})|\\
&\lesssim   \|Q_{0,b_1}^h\|_{\ell^2_{-1/2}} \big(\| Q_{a,b_2}\|_{\ell^2_{3/2}} \|Q_{0,b_3}\|_{\ell^2_{1+}}+\| Q_{a,b_2}\|_{\ell^2_{1+}} \|Q_{0,b_3}\|_{\ell^2_{3/2}}\big)\\
&\lesssim \|h\|_{\ell^2_{-1/2}} \, \eps^{b_1-} \big(\eps^{\frac{a}{2}+b_2-}\eps^{-\frac{1}{2}+b_3-} + \eps^{\frac12+\frac{a}{2}+b_2-}\eps^{-1+b_3-}\big)\\
&\lesssim \|h\|_{\ell^2_{-1/2}} \, \eps^{-\frac12+\frac{a}{2}+b-}.
\end{align*}
If $a_1\neq 0$, the worst case is when $a_1=a, a_2 = a_3 = 0$. We estimate
\begin{align*}
 |H_3(Q_{a,b_1}^h, Q_{0 ,b_2}, Q_{0,b_3})|
&\lesssim   \|Q_{a,b_1}^h\|_{\ell^{2}_{1/2-}}  \| Q_{0,b_2}\|_{\ell^2_{1}} \|Q_{0,b_3}\|_{\ell^2_{1/2+}} \\
&\lesssim \|h\|_{\ell^2_{-1/2}} \, \eps^{\frac{a}{2}+b_1-}  \eps^{-\frac12+b_2-}\eps^{ b_3-}  \\
&\lesssim \|h\|_{\ell^2_{-1/2}} \, \eps^{-\frac12+\frac{a}{2}+b-}.
\end{align*}

\end{proof}

Now we prove Propositions \ref{prop:H3}, \ref{prop:Qab}, \ref{prop:Qabh}.
\begin{proof}[Proof of Proposition \ref{prop:H3}] To prove the Proposition we will repeatedly use, without mentioning, the results of Lemma \ref{lem:calc}.
Since $H_3$ is symmetric in $q_1,q_2,q_3$, it suffices to consider the case $(j_1,j_2,j_3)=(1,2,3)$. We start with the second assertion:
\begin{align*}
|H_3(q_1,q_2,q_3)|&\lesssim   \sum_{n_1+n_2+n_3=0} \sqrt{n_1 n_2 n_3} \,\, |q_1(n_1)q_2(n_2)q_3(n_3)|
\\ &= \big\langle \sqrt{\cdot} \,|q_1|, \sqrt{\cdot} \, |q_2|  *\sqrt{\cdot}\,  |q_3| \big\rangle\lesssim \|\sqrt{\cdot} \, q_1\|_{\ell^{1+}} \|\sqrt{\cdot}\,  |q_2|  *\sqrt{\cdot}  \,|q_3|  \|_{\ell^{\infty-}}\\
& \lesssim \|  q_1\|_{\ell^2_{1- }} \|\sqrt{\cdot}\,  q_2 \|_{\ell^{2-}} \| \sqrt{\cdot}  \,q_3  \|_{\ell^{2-}} \lesssim \|  q_1\|_{\ell^2_{1- }} \|  q_2 \|_{\ell^{2 }_{1/2+}} \|  q_3  \|_{\ell^{2 }_{1/2+}}.
\end{align*}
To prove the other three assertions note that for $s<1/2$,
\begin{align} \nn
&|H_3(q_1,q_2,q_3)| \lesssim  \\ & \lesssim \sum_{n_1+n_2+n_3=0} |n_1|^s \big(|n_2|^{1 -s} |n_3|^{1/2}+|n_2|^{1/2 } |n_3|^{1 -s}\big) \,\,  |q_1(n_1)q_2(n_2)q_3(n_3)|\nonumber\\
&= \big\langle |\cdot|^s \,|q_1|, |\cdot|^{1-s}\, |q_2|  *|\cdot|^{1/2}\,  |q_3| \big\rangle+
\big\langle |\cdot|^s \,|q_1|, |\cdot|^{1/2 }\, |q_2|  *|\cdot|^{1-s}\,  |q_3| \big\rangle\nonumber
\end{align}
\be
\lesssim \|  q_1\|_{\ell^2_s} \Big( \big\| |\cdot|^{1-s}\, |q_2|  *|\cdot|^{1/2}\,  |q_3| \big\|_{\ell^2}
+\big\| |\cdot|^{1/2}\, |q_2|  *|\cdot|^{1-s}\,  |q_3| \big\|_{\ell^2}\Big). \label{eq:H3s}
\ee
For $s=0$, we bound \eqref{eq:H3s} by
$$
\|  q_1\|_{\ell^2 } \big(\| q_2  \|_{\ell^2_1}   \|q_3\|_{\ell^1_{1/2}}+ \| q_2  \|_{\ell^1_{1/2}}   \|q_3\|_{\ell^2_{1}}\big) \lesssim
 \|  q_1\|_{\ell^2 }  \| q_2  \|_{\ell^2_{1+}}   \|q_3\|_{\ell^2_{1+}},
$$
which proves the first assertion. For $s=-1/2$, we bound \eqref{eq:H3s}  by
\begin{multline}\nn
\|  q_1\|_{\ell^2_{-1/2} } \big(\| q_2  \|_{\ell^2_{3/2}}   \|q_3\|_{\ell^1_{1/2}}+ \| q_2  \|_{\ell^1_{1/2}}   \|q_3\|_{\ell^2_{3/2}}\big)
\\ \lesssim \|  q_1\|_{\ell^2_{-1/2} } \big(\| q_2  \|_{\ell^2_{3/2}}   \|q_3\|_{\ell^2_{1+}}+ \| q_2  \|_{\ell^2_{1+}}   \|q_3\|_{\ell^2_{3/2}}\big),
\end{multline}
which proves  the third assertion. Finally for $s=1/2-\delta$, $\delta>0$, we bound \eqref{eq:H3s}  by
\begin{multline}\nn
\|  q_1\|_{\ell^2_{1/2-\delta} } \big(\| q_2  \|_{\ell^{1+}_{1/2+\delta}}   \|q_3\|_{\ell^{2-}_{1/2}}+ \| q_2  \|_{\ell^{1+}_{1/2}}   \|q_3\|_{\ell^{2-}_{1/2+\delta}}\big)\\
\\ \lesssim \|  q_1\|_{\ell^2_{1/2-\delta} } \big(\| q_2  \|_{\ell^2_{1}}   \|q_3\|_{\ell^2_{1/2+}}+ \| q_2  \|_{\ell^2_{1 }}   \|q_3\|_{\ell^2_{1/2+\delta+}}\big).
\end{multline}
\end{proof}
\begin{proof}[Proof of Proposition \ref{prop:Qab}]
We start with the case $a=0$. Note that for $b=0$, the statement is true for any $s\in[0,3/2]$ since $Q_{0,0}=q$. For $b\geq 1$, we use a simple induction. We can write
$$Q_{0,b}=f_2(Q_{0,b_1},Q_{0,b_2},Q_{0,b_3}),$$
with $b_1+b_2+b_3=b-1$. Using Lemma~\ref{lem:f2}, for any $s\in [0,3/2]$, we have
$$\|Q_{0,b}\|_{\ell^2_s}\lesssim\sum \|Q_{0,b_{j_1}}\|_{\ell^2_s} \|Q_{0,b_{j_2}}\|_{\ell^2_{0+}} \|Q_{0,b_{j_3}}\|_{\ell^2_{0+}}, $$
where the sum is over all permutations $(j_1,j_2,j_3)$ of $(1,2,3)$. By the induction hypothesis the last sum can be estimated by
$$
\|Q_{0,b}\|_{\ell^2_s}\lesssim \sum \eps^{\frac12-s+b_{j_1}-}\eps^{\frac12 +b_{j_2}-}\eps^{\frac12 +b_{j_3}-}\lesssim \eps^{\frac12-s+b-}.
$$
In the case $a\geq 1$, we set up an induction on $a$. We first prove that the statement is valid for $a=1$ and for any $s\in [1,3/2]$, $b\geq 0$. We write
$$Q_{1,b}=f_1(Q_{0,b_1},Q_{0,b_2}),$$
$b_1+b_2=b$. By Lemma~\ref{lem:f1} we estimate
\begin{align*}
\|Q_{1,b}\|_{\ell^2_{3/2}}&=\|f_1(Q_{0,b_1},Q_{0,b_2})\|_{\ell^2_{3/2}}  \lesssim \|Q_{0,b_1}\|_{\ell^2_{1/2}} \|Q_{0,b_2}\|_{\ell^2_{0+}} + \|Q_{0,b_2}\|_{\ell^2_{1/2}} \|Q_{0,b_1}\|_{\ell^2_{0+}} \\
&\lesssim \eps^{b_1-}\eps^{\frac12+b_2-}+\eps^{b_2-}\eps^{\frac12+b_1-} \lesssim \eps^{\frac12+b-}.
\end{align*}
Again by Lemma~\ref{lem:f1} we estimate
\begin{align*}
\|Q_{1,b}\|_{\ell^2_{1-}}&=\|f_1(Q_{0,b_1},Q_{0,b_2})\|_{\ell^2_{1-}}  \lesssim \|Q_{0,b_1}\|_{\ell^2} \|Q_{0,b_2}\|_{\ell^2}  \lesssim \eps^{\frac12+b_1-}\eps^{\frac12+b_2-}  \lesssim \eps^{1+b-}.
\end{align*}
Now by a simple interpolation, for $s\in [1,3/2)$,
$$
\|Q_{1,b}\|_{\ell^2_{s}}\lesssim \|Q_{1,b}\|_{\ell^2_{1-}}^\theta \|Q_{1,b}\|_{\ell^2_{3/2}}^{1-\theta}
\lesssim \eps^{\theta(1+b-)}\eps^{(1-\theta)(\frac12+b-)}=  \eps^{\frac{1+\theta}{2}+b-},
$$
where $\theta=3-2s-$. This implies
$$
\|Q_{1,b}\|_{\ell^2_{s}} \lesssim \eps^{2-s+b-}.
$$
We proceed by induction on $a>1$. We have
$$Q_{a,b}=f_1(Q_{a_1,b_1},Q_{a_2,b_2}),$$
$a_1+a_2=a-1$, $b_1+b_2=b$. The worst case (in terms of gain in powers of $\eps$) is when $a_1=a-1$ and  $a_2=0$.
As above, using the induction hypothesis and Lemma~\ref{lem:f1} we have
\begin{align*}
\|Q_{a,b}\|_{\ell^2_{3/2}}&   \lesssim \|Q_{a-1,b_1}\|_{\ell^2_{1/2}} \|Q_{0,b_2}\|_{\ell^2_{0+}} + \|Q_{0,b_2}\|_{\ell^2_{1/2}} \|Q_{a-1,b_1}\|_{\ell^2_{0+}} \\
&\lesssim\|Q_{a-1,b_1}\|_{\ell^2_{1}} \|Q_{0,b_2}\|_{\ell^2_{1/2}} \lesssim \eps^{\frac{a}{2}+b_1-}\eps^{b_2-} = \eps^{\frac{a}{2}+b-}.
\end{align*}
Similarly,
\begin{align*}
\|Q_{a,b}\|_{\ell^2_{1-}}&  \lesssim \|Q_{a-1,b_1}\|_{\ell^2} \|Q_{0,b_2}\|_{\ell^2} \lesssim \|Q_{a-1,b_1}\|_{\ell^2_1} \|Q_{0,b_2}\|_{\ell^2} \\ &  \lesssim \eps^{\frac{a}{2}+b_1-}\eps^{\frac12+b_2-}  \lesssim \eps^{\frac12+\frac{a}{2}+b-}.
\end{align*}
The statement for $s\in [1,3/2)$ follows from interpolation as above.
\end{proof}

\begin{proof}[Proof of Proposition \ref{prop:Qabh}]
We give a proof only for the case $h\in\ell^2_{-1/2}$. The proof for the case $h\in \ell^2$ is essentially the same.  We start with the case $a=0$. Note that for $b=0$, the statement is true since $Q_{0,0}^h=h$.
For $b\geq 1$, we use a simple induction. We can write without loss of generality
$$Q_{0,b}^h=f_2(Q_{0,b_1}^h,Q_{0,b_2},Q_{0,b_3}),$$
with $b_1+b_2+b_3=b-1$. Using Lemma~\ref{lem:f2},  we have
\begin{align*}
\|Q_{0,b}^h\|_{\ell^2_{-1/2}}&\lesssim  \|Q_{0,b_{1}}^h\|_{\ell^2_{-1/2}} \|Q_{0,b_2}\|_{\ell^2_{0+}} \|Q_{0,b_3}\|_{\ell^2_{0+}}\\
&\lesssim \|h\|_{\ell^2_{-1/2}} \eps^{b_1-} \eps^{\frac12+b_2-} \eps^{\frac12+b_3-}
= \|h\|_{\ell^2_{-1/2}} \eps^{b-}.
\end{align*}
The second inequality follows from the induction hypothesis and Proposition~\ref{prop:Qab}.

In the case $a\geq 1$, we set up an induction on $a$. We first prove that the statement is valid for $a=1$  for any $b\geq 0$. We write, without loss of generality,
$$Q_{1,b}^h=f_1(Q_{0,b_1}^h,Q_{0,b_2}),$$
$b_1+b_2=b$. By Lemma~\ref{lem:f1} we estimate
\begin{align*}
\|Q_{1,b}^h\|_{\ell^2_{1/2-}}&=\|f_1(Q_{0,b_1}^h,Q_{0,b_2})\|_{\ell^2_{1/2-}}  \lesssim \|Q_{0,b_1}^h\|_{\ell^2_{-1/2}} \|Q_{0,b_2}\|_{\ell^2}  \\
&\lesssim \|h\|_{\ell^2_{-1/2}} \eps^{b_1-}\eps^{\frac12+b_2-}  \lesssim \|h\|_{\ell^2_{-1/2}} \eps^{\frac12+b-}.
\end{align*}
The second inequality follows from the first part of the proof and Proposition~\ref{prop:Qab}.

We proceed by induction on $a>1$. We have, without loss of generality,
$$Q_{a,b}^h=f_1(Q_{a_1,b_1}^h,Q_{a_2,b_2}),$$
$a_1+a_2=a-1$, $b_1+b_2=b$.
Using Lemma~\ref{lem:f1},  we have
\begin{align*}
\|Q_{a,b}^h\|_{\ell^2_{1/2-}}&  \lesssim \|Q_{a_1,b_1}^h\|_{\ell^2_{-1/2}} \|Q_{a_2,b_2}\|_{\ell^2}  \\
&\lesssim \|h\|_{\ell^2_{-1/2}} \eps^{\frac{a_1}{2}+b_1-}\eps^{\frac12+\frac{a_2}{2}+b_2-}  \lesssim \|h\|_{\ell^2_{-1/2}} \eps^{\frac{a}{2}+b-}.
\end{align*}
The second inequality follows from the induction hypothesis and Proposition~\ref{prop:Qab} by considering the cases $a_1=0$, $a_1\neq 0$ and  $a_2=0$, $a_2\neq 0$.
\end{proof}

\subsection{Remainder estimates} \label{sec:rmn}
In this section we estimate the error terms involving integrals. By \eqref{eq:R1def} and \eqref{eq:Rfinal}, it suffices to prove the following inequalities
\ba \label{eq:R_int1}
&\sup_{\tau\in [0,1]} \big\|\partial_q \,(g_{F_1}^3 H_3 \circ \Phi_{F_1}^\tau) \big\|_{\ell^2_s}\lesssim \eps^{1-s-},\quad\quad \,\,s\in [0,1/2]\\
& \,\,\,\,\,\,\,\,  \sup_{\tau\in [0,1]} \big\|\partial_q \,g^{}_{F_2}(g_{F_1}^3 H_3 \circ \Phi_{F_1}^\tau)   \big\|_{\ell^2_s}\lesssim \eps^{1-s-},\quad\quad \,\,s\in [0,1/2]\label{eq:R_int2}\\
&\sup_{\tau\in [0,1]} \big\| \partial_q\,\big( g_{F_2}^2 (H\circ \Phi_{F_1}^1) \circ \Phi_{F_2}^\tau \big) \big\|_{\ell^2_s}\lesssim \eps^{1-s-},\quad\quad s\in [0,1/2].\label{eq:R_int3}
\ea
To prove \eqref{eq:R_int1}, let $w=\Phi_{F_1}^\tau(q)$.
Note that
\begin{align*}
\big\| \partial_q   \big( g_{F_1}^a H_3 \circ \Phi_{F_1}^\tau\big) \big\|_{\ell^2_s}   & = \sup_{\|h\|_{\ell^2_{-s}}=1} \Big|\sum_{m,k} \frac{\partial g_{F_1}^a H_3 }{\partial w(m) } \frac{\partial w(m)}{\partial q(-k) } \,h(k) \Big| \\
& = \sup_{\|h\|_{\ell^2_{-s}}=1} \Big|\sum_{m} \frac{\partial g_{F_1}^a H_3 }{\partial w(m) } \big(\sum_k\frac{\partial w(m)}{\partial q(-k) } \,h(k) \big)\Big|\\
&\lesssim \big\|\partial_w g_{F_1}^a H_3 \big\|_{\ell^2_s} \sup_{\|h\|_{\ell^2_{-s}}=1} \big\|\sum_{k}   \frac{\partial w(m)}{\partial q(-k) } \,h(k) \big\|_{\ell^2_{-s}}.
\end{align*}
Since $F_1$ is near identity, by our assumptions on $q$, $\|w\|_{\ell^2_s}\lesssim \eps^{\frac12-s}$ for $s\in[0,3/2]$. Therefore, Theorem~\ref{thm:H3ab} implies that $\|\partial_w g_{F_1}^a H_3 \big\|_{\ell^2_s} \lesssim \eps^{1-s-}$ (for $a\geq 2$). Thus, it suffices to prove that
\ba\label{eq:Dkw}
 \|T(h)  \|_{\ell^2_{-s}} \lesssim \|h\|_{\ell^2_{-s}},\,\,\,\,\,\,s\in [0,1/2]
\ea
where $T(h):=\sum_{k} D_k w(m)  \,h(k)$, and    $  D_k w(m):=\frac{\partial w(m)}{\partial q(-k) } $.

To prove \eqref{eq:Dkw} first note that $w(m)$ is the solution at $t=\tau$ of the system
$$
\frac{d w(m)}{dt}= \frac{\partial F_1}{\partial w(-m)}=f_1(w,w)(m),\quad \quad w|_{t=0}=q.
$$
Differentiating this equation with respect to initial condition $q(-k)$, we see that
$$
\frac{d \, D_kw(m)}{dt}= 2f_1(D_kw,w),\quad \quad D_k w(m)|_{t=0}=\delta_{-k,m}.
$$
Pairing both sides with $h(k)$, we have the following equation for $T(h)$
$$
\frac{d \, T(h)(m)}{dt}= 2f_1(T(h),w),\quad \quad T(h)(m)|_{t=0}=h(-m).
$$
Therefore \eqref{eq:Dkw} is satisfied at $\tau=0$, and by Gronwall's lemma (since $\tau\in[0,1]$), it suffices to see that
$\|f_1(T(h),w)\|_{\ell^2_{-s}}\lesssim \|T(h)\|_{\ell^2_{-s}}$ for $s\in[0,1/2]$. This immediately follows from Lemma~\ref{lem:f1}.

The remaining estimates \eqref{eq:R_int2} and \eqref{eq:R_int3} follow from similar considerations using
Lemma~\ref{lem:f1}, Lemma~\ref{lem:f2}, and Theorem~\ref{thm:H3ab}. We omit the details.


\begin{thebibliography}{100}
\bibitem{Bou1} J.~Bourgain,{\em Fourier transform restriction phenomena for certain lattice subsets and
applications to nonlinear evolution equations. Part I: Schr\"odinger equations}, GAFA, {\bf 3} No. 2 (1993),
107--156.
\bibitem{Bou2} J.~Bourgain, {\em Fourier transform restriction phenomena for certain lattice subsets and
applications to nonlinear evolution equations. Part II: The KdV equation}, GAFA, {\bf 3} (1993), 209--262.
\bibitem{ckstt} J.~Colliander, M.~Keel, G.~Staffilani, H.~Takaoka, T.~Tao, {\it Weakly turbulent solutions for the cubic defocusing nonlinear Schrödinger equation,} preprint, http://arxiv.org/abs/0808.1742.
\bibitem{erdogan_vz} M.~B.~Erdo\u{g}an, V.~Zharnitsky, {\em Quasi-linear dynamics in nonlinear Schr\" odinger equation with periodic boundary conditions}, Commun.\ Math.\ Phys.\,  {\bf 281} (2008),  655--673.
\bibitem{Naum} N.~Hayashi, P.~Naumkin, {\em Asymptotics for large time of solutions to the nonlinear Schr\"odinger and Hartree equations}, Amer.\ Jour.\ Math.\, {\bf 120} (1998), 369--389.
\bibitem{KappelerPoschel} T.~Kapeller, J.~P\"oschel, {\em KdV and KAM},   A Series of Modern
              Surveys in Mathematics   {\bf 45}, Springer-Verlag, Berlin, (2003).
\bibitem{KharifPelin} C. Kharif, E. Pelinovsky, {\em Physical mechanisms of rogue wave phenomenon}, Eur.\ Jour.\ Mech.\ B/Fluids, {\bf 22} (2003) 603--634.
\bibitem{Kuk} S.~Kuksin, {\em Analysis of Hamiltonian PDEs}, Oxford University Press, New York, 2000.
\bibitem{Temam} R.~Temam, {\em Sur un probleme non lineaire}, J.\ Math.\ Pures Appl.\ {\bf 48} (1969), 159--172.
\end{thebibliography}
\end{document}